\documentclass[12pt]{amsart}
\usepackage[cp1251]{inputenc}
\usepackage[english]{babel}
\usepackage{amsmath}
\usepackage{amssymb}
\usepackage{amsfonts}

\begin{document}


\title[On some solvability theorems]
      {On some solvability theorems for pseudo-differential equations}
\author{Vladimir Vasilyev}
\address{Chair of Applied Mathematics and Computer Modeling\\
    Belgorod State National Research University\\
         Pobedy street 85, Belgorod 308015, Russia}

        \email{vbv57@inbox.ru}

  \author{Victor Polunin}
\address{Chair of Applied Mathematics and Computer Modeling\\
    Belgorod State National Research University\\
         Pobedy street 85, Belgorod 308015, Russia}

        \email{polunin@bsu.edu.ru}

 \author{Igor Shmal}
\address{Chair of Applied Mathematics and Computer Modeling\\
    Belgorod State National Research University\\
         Pobedy street 85, Belgorod 308015, Russia}

\email{1247978@bsu.edu.ru}

\begin{abstract}
We study a model elliptic pseudo-differential
equation and simplest boundary value problems for
a half-space and a special cone in
Sobolev--Slobodetskii spaces which have
different smoothness with respect to separate
variables. Sufficient conditions for a unique
solvability for such boundary value problems are
described.
\end{abstract}

\maketitle

\section{Introduction}

The theory of pseudo-differential operators was appeared near a half-century ago, and it has taken attention of mathematicians for a long time \cite{TRE,TAY,HER}. More general Fourier integral operators and new functional spaces were studied in this context. As a rule the theory means constructing a symbolic calculus and an index formula. Such a theory is very convenient for generalization on smooth compact manifolds without a boundary, but for more complicated situations new constructions and new approaches were needed. More complicated situations mean presence of a smooth boundary, or more generally a non-smooth boundary. For manifolds  with a smooth boundary a certain approach was suggested in \cite{E}, and it was based the factorization principle for an elliptic symbol at boundary points. This method is not applicable for manifolds with a non-smooth boundary, and it has initiated a lot of approaches for "non-smooth" situations \cite{SCH1,SCH2,NIS,RUZ,PLA,RAB,V0}.

This paper presents a future development of the second author's approach \cite{V3,V4,V1,V2} to the theory of pseudo-differential equations and related boundary value problems in non-smooth domains.

\section{Elliptic pseudo-differential operators}

\subsection{Sobolev--Slobodetskii spaces
of different smoothness}

Following to \cite{NRSW} (see also \cite{VG}) we introduce useful notations. A multidimensional Euclidean space
${\mathbb R}^M$ is represented as an orthogonal sum of subspaces in which only some of coordinates  $x_1,x_2,\dots,x_M$ are nor vanishing. Namely, if
 $K\subset{1,...,M}$ is not empty set we put
\[
{\mathbb R}^K=\{x\in{\mathbb R}^M: x=(x_1,\dots,x_M),
x_j=0, \forall j\notin K\}\subset{\mathbb R}^M.
\]

Let $K_1,K_2,\dots,K_n\subset\{1,2,\dots,M\}$ be a nonempty set so that
\[
\bigcup\limits_{j=1}^nK_j=\{1,2,\dots,M\},~~~K_i\cap
K_j=\emptyset, i\neq j, ~~~card~ K_j=k_j.
\]

Thus, we obtain the representation
\[
{\mathbb R}^M={\mathbb R}^{K_1}\oplus{\mathbb R}^{K_2}\oplus\dots\oplus{\mathbb R}^{K_n},
\]
where $x_{K_j}$ is an element of the space ${\mathbb R}^{K_j}$. For functions defined in ${\mathbb R}^M$ we use the standard Fourier transform
\[
\tilde u(\xi)=\int\limits_{{\mathbb R}^M}e^{ix\cdot\xi}u(x)dx,~~~\xi=(\xi_1,\dots,\xi_M).
\]

Let $S=(s_1,\dots,s_n)$. Now we introduce the Sobolev--Slobodetskii space
$H^{S}({\mathbb R}^M)$ as a Hilbert space with the inner product
\[
(f,g)=\int\limits_{{\mathbb R}^M}f(x)\overline{g(x)}dx
\]
and the norm
\[
||f||_S=\left(\int\limits_{{\mathbb R}^M}(1+|\xi_{K_1}|)^{2s_1}(1+|\xi_{K_2}|)^{2s_2}\cdots(1+|\xi_{K_n}|)^{2s_n}|\tilde
f(\xi)|^2d\xi\right)^{1/2}.
\]

Such $H^S$-spaces have the same properties similar to usual Sobolev--Slobodetskii spaces \cite{VG}. Particularly, the usual space $H^s({\mathbb R}^M)$ is obtained under the following choice of subsets  $K_j$ and parameters $s_j$:
\[
K_1=K_2=\dots=K_{n-1}=\emptyset,~~~K_n=\{1,2,\dots,M\},~~~S=(0,0,\dots,0,s).
\]

\subsection{Model operators and equations}

According to the local principle we will concentrate on studying a model pseudo-differential equation with operator with a symbol non-depending on a spatial variable.

\subsubsection{Model pseudo-differential operators}

Let $\tilde A(\xi), \xi\in{\mathbb R}^M$ be a measurable function. A model pseudo-differential operator $A$ is defined as follows
\[
(Au)(x0=\frac{1}{(2\pi)^M}\int\limits_{{\mathbb R}^M}\int\limits_{{\mathbb R}^M}e^{i(x-y)\cdot\xi}\tilde A(\xi)u(y)dyd\xi,
\]
and the function $\tilde A(\xi)$ is called a symbol of the pseudo-differential operator $A$.

We consider here the following class of symbols $A(\xi)$ satisfying the condition
$$
\begin{array}{rcl}
c_1\prod\limits_{j=1}^n(1+|\xi_{K_j}|)^{\alpha_j}\leq|A(\xi)|\leq c_2\prod\limits_{j=1}^n(1+|\xi_{K_j}|)^{\alpha_j},\\
~~~\alpha_j\in{\mathbb R}, j=1,2,\dots,n,
\end{array}
\eqno(*)
$$
with positive constants $c_1,c_2$.

Let $\alpha=(\alpha_1,\dots,\alpha_n)$.

{\bf Lemma 1.} {\it
Let $A$ be a pseudo-differential operator with the symbol $\tilde A(\xi)$ satisfying the condition $(*)$. Then $A: H^S({\mathbb R}^M)\rightarrow H^{S-\alpha}({\mathbb R}^M)$ is a linear bounded operator.
}

\begin{proof}
Indeed, we have
\[
||Au||^2_{S-\alpha}=\int\limits_{{\mathbb R}^M}\prod\limits_{j=1}^n(1+|\xi_{K_j}|)^{2(s_j-\alpha_j)}\widetilde
{Au}(\xi)|^2d\xi=
\]
\[
=\int\limits_{{\mathbb R}^M}(1+|\xi_{K_1}|)^{2(s_1-\alpha_1)}(1+|\xi_{K_2}|)^{2(s_2-\alpha_2)}\cdots(1+|\xi_{K_n}|)^{2(s_n-\alpha_n)}|\tilde
A(\xi)\tilde u(\xi)|^2d\xi\leq
\]
\[
\leq c_2\int\limits_{{\mathbb R}^M}(1+|\xi_{K_1}|)^{2s_1}(1+|\xi_{K_2}|)^{2s_2}\cdots(1+|\xi_{K_n}|)^{2s_n}|\tilde
u(\xi)|^2d\xi=c_2||u||^2_S,
\]
and the proof is completed.
\qed
\end{proof}

Thus, we can start studying a solvability for the equation
\begin{equation}\label{1}
(Au)(x)=v(x),~~~x\in{\mathbb R}^M,
\end{equation}
where  $A$ is a pseudo-differential operator with the symbol $\tilde A(\xi)$ satisfying the condition $(*)$, and the right hand side $v\in H^{S-\alpha}({\mathbb R}^M)$.

{\bf Corollary 1.} {\it
If $A$ is a pseudo-differential operator with the symbol $\tilde A(\xi)$ satisfying the condition $(*)$ then the equation $(\ref{1})$ with an arbitrary right hand side $v\in H^{S-\alpha}({\mathbb R}^M)$ has unique solution $u\in H^S({\mathbb R}^M)$. The a priori estimate
\[
||u||_S\leq C||v||_{S-\alpha}
\]
holds.
}

\begin{proof}
The operator $A^{-1}$ with the symbol $\tilde A^{-1}(\xi)$ is a pseudo-differential operator. Its symbol satisfies the condition $(*)$ with order $-\alpha$ instead of $\alpha$. Then we have
\[
\tilde u=\tilde A^{-1}\tilde v,
\]
and therefore
\[
||u||^2_S=||A^{-1}v||^2_S=
\]
\[
\int\limits_{{\mathbb R}^M}(1+|\xi_{K_1}|)^{2s_1}(1+|\xi_{K_2}|)^{2s_2}\cdots(1+|\xi_{K_n}|)^{2s_n}|\tilde A^{-1}(\xi)\tilde
v(\xi)|^2d\xi\leq c_1^{-2}||v||^2_{S-\alpha},
\]
and the sentence is proved.
\qed
\end{proof}

Unfortunately, such a simple conclusion is possible for the space ${\mathbb R}^M$. If we will take a domain $D\subset{\mathbb R}^M$ and will try to study a solvability for similar equation then we will obtain a lot of difficulties related to invertibility of operators.

We extract special {\it canonical} domains $D$ in Euclidean space ${\mathbb R}^M$. Such domains are conical domains and we will start from a standard convex cone in Euclidean space non-including a whole straight line. Let $C_{K_j}\subset {\mathbb R}^{K_j}$ and we would like to consider the equation
\begin{equation}\label{2}
(Au)(x)=v(x),~~~x\in C_{K_j}.
\end{equation}
Direct applying the Fourier transform does not give the required answer since we have no the convolution theorem. The equation (\ref{2}) can be rewritten in the form
\[
(P_{K_j}Au)(x)=v(x),~~~x\in  C_{K_j},
\]
where $P_{K_j}$ is the restriction on $C_{K_j}$,
$$
(P_{K_j}u)(x)=\left\{
\begin{array}{rcl}
u(x), &x\in  C_{K_j};\\
0, &x\notin\overline{ C}_{K_j}.
\end{array}
\right.
$$
and to use the Fourier transform we need to know what is the Fourier image of the operator $P_{K_j}$.

\subsubsection{Structure of projectors}

For a general convex cone $C^m\subset{\mathbb R}^m$ one can define the Bochner kernel \cite{BM,VL1,VL2}
\[
B_m(z)=\int\limits_Ce^{ix\cdot z}dx,~~~z=(z_1,\dots,z_m),
\]
and the following representation in Fourier imaged
\[
(FP_+u)(\xi)=\lim\limits_{\tau\to 0+}\int\limits_{{\mathbb R}^m}B_m(\xi'-\eta',\xi_m-\eta_m+i\tau)\tilde u(\eta',\eta_m)d\eta,
\]
here $P_+$ is the projector on the cone $C^m$ \cite{V0,V1}. There are certain concrete realizations in the latter formula.

{\bf Example 1.}

We consider here one-dimensional case in which we have only one cone, and this cone is ${\mathbb R}_+$ \cite{E}. For this case it was proved for a function $u(x), x\in{\mathbb R},$ that
\[
(FP_+u)(\xi)=\frac{1}{2}\tilde u(\xi)+\frac{i}{2\pi}p.v.\int\limits_{-\infty}^{+\infty}\frac{\tilde u(\eta)d\eta}{\xi-\eta}.
\]

As a consequence we have for a function $u(x),x\in{\mathbb R}^m$ and the cone ${\mathbb R}^m_+=\{x\in{\mathbb R}^m: x=(x',x_m), x_m>0\}$ the following result
\[
(FP_+u)(\xi)=\frac{1}{2}\tilde u(\xi)+\frac{i}{2\pi}p.v.\int\limits_{-\infty}^{+\infty}\frac{\tilde u(\xi',\eta_m)d\eta_m}{\xi_m-\eta_m},~~~\xi=(\xi',\xi_m).
\]

{\bf Example 2.}

Let $m=2$, and
\[
C^a_+=\{x\in{\mathbb R}^2: x=(x_1,x_2), x_2>a|x_1|, a>0\}.
\]
Then we have \cite{V5}
$$
(FP_{C^a_+}u)(\xi)=\frac{\tilde u(\xi_1+a\xi_2,\xi_2)+\tilde u(\xi_1-a\xi_2,\xi_2)}{2}+
$$
$$
+v.p.\frac{i}{2\pi}\int\limits_{-\infty}^{+\infty}\frac{\tilde u(\eta,\xi_2)d\eta}{\xi_1+a\xi_2-\eta}-
v.p.\frac{i}{2\pi}\int\limits_{-\infty}^{+\infty}\frac{\tilde u(\eta,\xi_2)d\eta}{\xi_1-a\xi_2-\eta}.
$$

{\bf Example 3.}
Let $m=3$, and   $C_+^{a_1a_2}=\{x\in{\mathbb R}^3: x_2>a_1|x_1|+a_2|x_2|\}$. Then
\[
(FP_{C_+^{a_1a_2}}u)(\xi_1,\xi_2,\xi_3)=
\]
\[
=\frac{\tilde u(\xi_1-a_1\xi_3,\xi_2-a_2\xi_3,\xi_3)+\tilde u(\xi_1+a_1\xi_3,\xi_2-a_2\xi_3,\xi_3)}{4}+
\]
\[
+\frac{1}{2}(S_1\tilde u)(\xi_1+a_1\xi_3,\xi_2-a_2\xi_3,\xi_3)-\frac{1}{2}(S_1\tilde u)(\xi_1-a_1\xi_3,\xi_2-a_2\xi_3,\xi_3)+
\]
\[
+\frac{\tilde u(\xi_1-a_1\xi_3,\xi_2+a_2\xi_3,\xi_3)+\tilde u(\xi_1+a_1\xi_3,\xi_2+a_2\xi_3,\xi_3)}{4}+
\]
\[
+\frac{1}{2}(S_1\tilde u)(\xi_1+a_1\xi_3,\xi_2+a_2\xi_3,\xi_3)-\frac{1}{2}(S_1\tilde u)(\xi_1-a_1\xi_3,\xi_2+a_2\xi_3,\xi_3)+
\]
\[
+\frac{(S_2\tilde u)(\xi_1-a_1\xi_3,\xi_2+a_2\xi_3,\xi_3)+(S_2\tilde u)(\xi_1+a_1\xi_3,\xi_2+a_2\xi_3,\xi_3)}{2}+
\]
\[
+(S_1S_2\tilde u)(\xi_1+a_1\xi_3,\xi_2+a_2\xi_3,\xi_3)-(S_1S_2\tilde u)(\xi_1-a_1\xi_3,\xi_2+a_2\xi_3,\xi_3)-
\]
\[
-\frac{(S_2\tilde u)(\xi_1-a_1\xi_3,\xi_2-a_2\xi_3,\xi_3)-(S_2\tilde u)(\xi_1+a_1\xi_3,\xi_2-a_2\xi_3,\xi_3)}{2}-
\]
\[
-(S_1S_2\tilde u)(\xi_1+a_1\xi_3,\xi_2-a_2\xi_3,\xi_3)+(S_1S_2\tilde u)(\xi_1-a_1\xi_3,\xi_2-a_2\xi_3,\xi_3).
\]
where
\[
(S_1u)(\xi_1,\xi_2,\xi_3)=v.p\frac{i}{2\pi}\int\limits_{-\infty}^{+\infty}\frac{u(\tau,\xi_2,\xi_3)d\tau}{\xi_1-\tau},
\]
\[
(S_2u)(\xi_1,\xi_2,\xi_3)=v.p\frac{i}{2\pi}\int\limits_{-\infty}^{+\infty}\frac{u(\xi_1,\eta,\xi_3)d\eta}{\xi_2-\eta}.
\]

 This case was studied in \cite{V6}

\subsubsection{Elliptic equations and complex variables}

This approach is related to the function theory of many complex variables, namely to functions which are holomorphic in radial tube domains \cite{BM,VL1,VL2}.

Let $C_{K_j}\subset{\mathbb R}^{K_j}, j=1,\dots,n,$ be convex cones non-including a whole straight line in ${\mathbb R}^{K_j}$. Let us compose the set $C=C_{K_1}\times\dots C_{K_n}$.

{\bf Lemma 2.} {\it
The set $C$ is a cone in ${\mathbb R}^M$ non-including a whole straight line in ${\mathbb R}^M$.
}

\begin{proof}
Indeed, $C$ is a cone since each $C_j$ is a cone. If we will assume that $C$ includes a certain line in ${\mathbb R}^M$ then we will conclude that each cone $C_j$ includes a certain straight line.
\qed
\end{proof}

Now we will start studying a solvability of the equation
\begin{equation}\label{3}
(Au)(x)=v(x),~~~x\in C,
\end{equation}
and the solution is sought in the space $H^S(C)$.
{\bf Definition 1.} {\it
The space $H^S(C)$ consists of functions (distributions) from\\
 $H^S({\mathbb R}^M)$ with supports in $\overline{C}$.
}
The right-hand side $v$ is chosen from the space $H^{S-\alpha}_0(C)$; by definition the space $H^S_0(C)$ is a space of distributions on
 $C$, admitting a continuation on $H^S({\mathbb R}^M)$. The norm in the space  $H^S_0(C)$ is defined as
$$
||v||^+_S=\inf ||\ell f||_S,
$$
where the {\it infimum} is taken over all continuations $\ell lf$ on the whole ${\mathbb R}^M$.

Fourier image of the the space $H^S(C)$ will be denoted by $\tilde H^S(C)$

{\bf Definition 2.} {\it
A radial tube domain over the cone   $C$ is called a domain in $M$-dimensional complex space ${\bf C}^M$ of the following type
$$
T(C)\equiv\{z\in{\bf C}^M: z=x+iy, x\in{\mathbb R}^M, y\in C\}.
$$

A conjugate cone  $\stackrel{*}{C}$ is called such a cone in which for all points the condition
$$
x\cdot y>0,~~~\forall y\in C,
$$
holds;  $x\cdot y$ means inner product for   $x$ and $y$.
}

{\bf Definition 3.} {\it
The wave factorization of an elliptic symbol $A(\xi)$ with respect to the cone $C$ is called its representation in the form
$$
A(\xi)=A_{\neq}(\xi)A_=(\xi),
$$
where factors  $A_{\neq}(\xi),A_=(\xi)$ must satisfy the following conditions:

1) $A_{\neq}(\xi),A_=(\xi)$ are defined for all$\xi\in{\mathbb R}^M$ may be except  the points $\xi\in\partial\stackrel{*}{C}$;

2) $A_{\neq}(\xi),A_=(\xi)$ admit an analytic continuation into radial tube domains $T(\stackrel{*} {C}),T(-\stackrel{*} {C})$ respectively with estimates
$$
|A_{\neq}^{\pm 1}(\xi+i\tau)|\leq c_1\prod\limits_{j=1}^n(1+|\xi_{K_j}|+|\tau_{K_j}|)^{\pm\ae_j},
$$
$$
|A_{=}^{\pm 1}(\xi-i\tau)|\leq c_2\prod\limits_{j=1}^n(1+|\xi_{K_j}|+|\tau_{K_j}|)^{\pm(\alpha_j-\ae_j)},~\forall\tau\in\stackrel{*} {C},~~~\ae_j\in{\mathbb R}.
$$

The vector $\ae=(\ae_1,\dots,\ae_n)$ is called an index of the wave factorization.
}

To apply the Fourier transform to the equation (\ref{3}) we need to know what is $FP_C$; here $F$ denotes the Fourier transform in $M$-dimensional space. Let us introduce the following notations. For every $C_{K_j}$ we consider corresponding radial tube domain $T(\stackrel{*}{C}_{K_j})$ over the conjugate cone and an element of $T(\stackrel{*}{C}_{K_j})$ will be denoted by $\xi_{K_j}+i\tau_{K_j}$. Moreover, for $\xi_{K_j}$ we will use the notation $\xi_{K_j}=(\xi'_{K_j},\xi_{k_j})$, where $\xi_{k_j})$ is the $k_j$th coordinate, and $\xi'_{K_j}$ denotes left other coordinates. The same notations will be used for $x\in{\mathbb R}^{K_j}, x_{K_j}=(x'_{K_j},x_{k_j})$.

As before we denote by $P_C$ the restriction operator on $C$. Obviously,
\[
P_C=\prod\limits_{j=1}^nP_{K_j}.
\]
and then
\[
B_M(z)=\prod\limits_{j=1}^nB_{k_j}(z_{K_j}),~~~z=(z_{K_1},\dots,z_{K_n}).
\]

The last our observation is the following:
\[
T(\stackrel{*}{C})=\prod\limits_{j=1}^nT(\stackrel{*}{C}_{K_j}),
\]
and the Bochner kernel $B_M(z)$ will be a holomorphic function in $T(\stackrel{*}{C})$.

{\bf Theorem 1.} {\it
If the symbol $A(\xi)$ admits the wave factorization with respect to the cone $C$ with the index $\ae$ such that $|\ae_j-s_j|<1/2, j=1,\dots,n,$ then the equation (\ref{3}) has unique solution in the space $H^S(C)$ for arbitrary right hand side $v\in H_0^{S-\alpha}(C)$.

The a priori estimate
\[
||u||_S\leq const~||v||^+_{S-\alpha}
\]
holds.
}

\begin{proof}
We use the Wiener--Hopf method \cite{E,V0}. Let $\ell v$ be an arbitrary continuation of $v$ onto ${\mathbb R}^M$. then we put
\[
u_-(x)=(\ell v)(x)-(Au)(x),
\]
so that $v_-(x)=0$ for $x\in C$. Further,
\[
(Au)(x)+u_-(x)=(\ell v)(x),
\]
and after applying the Fourier transform and the wave factorization we obtain
\begin{equation}\label{4}
A_{\neq}(\xi)\tilde u(\xi)+A^{-1}_{=}(\xi)\tilde u_-(\xi)=A^{-1}_{=}(\xi)\widetilde{(\ell v)}(\xi)
\end{equation}

Now we can use the following result (see \cite{V0}).
{\it Property 1.} {\it
If $\tilde u_-\in\widetilde H^S({\mathbb R}^M\setminus C), A^{-1}_=$ is a factor of the wave factorization then $A^{-1}_=\tilde u_-\in\widetilde H^{S+\alpha-\ae}({\mathbb R}^M\setminus C).$
}

Obviously, the summand $A_{\neq}(\xi)\tilde u(\xi)$ belongs to $\widetilde H^{S-\ae}(C)$ according to Lemma 1 and holomorphic properties, and $A^{-1}_{=}(\xi)\tilde u_-(\xi)$ belongs to $\widetilde H^{S-\ae}({\mathbb R}^M\setminus C)$ according to Property 1.

The right hand side $A^{-1}_{=}(\xi)\widetilde{(\ell v)}(\xi)$ belongs to the space $\widetilde H^{S-\ae}({\mathbb R}^M)$ (Lemma 1), and since  $|\ae_j-s_j|<1/2, j=1,\dots,n,$ it can be uniquely represented as
\begin{equation}\label{7}
 A^{-1}_{=}(\xi)\widetilde{(\ell v)}(\xi)=\tilde v_+(\xi)+\tilde v_-(\xi),
\end{equation}
 where
 \[
\tilde v_+(\xi)= B_M\left(A^{-1}_{=}(\xi)\widetilde{(\ell v)}(\xi)\right),~~~\tilde v_-(\xi)=(I-B_M)\left(A^{-1}_{=}(\xi)\widetilde{(\ell v)}(\xi)\right).
\]

The representation (\ref{7}) is true since the operator $B_M: \widetilde H^{\delta}({\mathbb R}^M)\rightarrow\widetilde H^{\delta}({\mathbb R}^M)$ for $|\delta_J|<1/2, j=1,\dots,n$, and we remind that $|\ae_j-s_j|<1/2, j=1,\dots,n,$.

Further, we rewrite the equality (\ref{4}) in the form
\[
A_{\neq}(\xi)\tilde u(\xi)-\tilde v_+(\xi)=\tilde v_-(\xi)-A^{-1}_{=}(\xi)\tilde u_-(\xi),
\]
and we obtain that a distribution from $H^{\delta}(C)$ equals to a distribution from $H^{\delta}({\mathbb R}^M\setminus\overline{C})$. But for such small $\delta$ this common distribution should be zero only \cite{V0}. Thus,
\[
A_{\neq}(\xi)\tilde u(\xi)-\tilde v_+(\xi)=0,
\]
or in other words
\[
\tilde u(\xi)= A^{-1}_{\neq}(\xi)B_M\left(A^{-1}_{=}(\xi)\widetilde{(\ell v)}(\xi)\right).
\]

A priori estimate is  based on Lemma 1 and boundedness property of the operator $B_M: \widetilde H^{\delta}({\mathbb R}^M)\rightarrow\widetilde H^{\delta}({\mathbb R}^M)$. Indeed,
\[
||u||_S=||\tilde u||_S=||A^{-1}_{\neq}(\xi)B_M\left(A^{-1}_{=}(\xi)\widetilde{(\ell v)}(\xi)\right)||_S\leq
\]
\[
\leq~const~||B_M\left(A^{-1}_{=}(\xi)\widetilde{(\ell v)}(\xi)\right)||_{S-\ae}\leq~const~||A^{-1}_{=}(\xi)\widetilde{(\ell v)}(\xi)||_{S-\ae}\leq
\]
\[
\leq~const~||\widetilde{(\ell v)}(\xi)||_{S-\alpha}=const~||\ell v||_{S-\alpha}\leq~const~||v||^+_{S-\alpha},
\]
and Theorem 1 is proved.
\qed
\end{proof}

\subsubsection{Multiply solutions}

For the cone $C_{K_j}, j=1,\dots,n,$  we suppose that a surface of this cone is given by the equation   $x_{k_j}=\varphi_j(x'_{K_j})$, where $\varphi_j:{\mathbb R}^{k_j-1}\rightarrow{\mathbb R}$ is a smooth function in  ${\mathbb R}^{k_j-1}\setminus\{0\}$, and $\varphi_j(0)=0$.

Let us introduce the following change of variables
$$
\left\{
\begin{array}{rcl}
t'_{K_j}&=&x'_{K_j}\\
t_{k_j}&=&x_{k_j}-\varphi_j(x'_{K_j})
\end{array}
\right.
$$
and we denote this operator by $T_{\varphi_j}:{\mathbb R}^{K_j}\rightarrow{\mathbb R}^{K_j}$. Since the cone is in one part of a half-space then points of the second part of a half-space will be fixed. Such change of variables can be defined for distributions also \cite{V6}.

Below we will use notation $F_m$ for the Fourier transform in $m$-dimensional space, so that the notation $F_{K_j}$ will be the Fourier transform in ${\mathbb R}^{K_j}$.

Following to \cite{V6} we conclude
\[
F_{K_j}T_{\varphi_j}=V_{\varphi_j}F_{K_j}.
\]
Further, we introduce $T_{\varphi}: {\mathbb R}^M\rightarrow{\mathbb R}^M$ by the formula
\[
T_{\varphi}=\prod\limits_{j=1}^nT_{\varphi_j}
\]
and construct the operator
\[
V_{\varphi}=\prod\limits_{j=1}^nV_{\varphi_j},
\]
for which we have
\[
F_{M}T_{\varphi}=V_{\varphi}F_{M}.
\]

Let us introduce vectors $N=(n_1,\dots,n_n), L=(l_1,\dots,l_n), \delta=(\delta_1,\dots,\delta_n), \\
n_j, l_j\in{\mathbb N}, |\delta_j|<1/2, j=1,\dots,n,$ and a polynomial $Q_N(\xi), \xi\in{\mathbb R}^M$ satisfying the condition
\begin{equation}\label{5}
|Q_N(\xi)|\sim\prod\limits_{j=1}^n(1+|\xi_{K_j}|)^{n_j},
\end{equation}

{\bf Theorem 2.} {\it
 If the symbol $A(\xi)$ admits the wave factorization with the index $\ae, \ae-S=N+\delta, $  then a  general solution of the equation  (\ref{3}) in Fourier images is given by the formula
\[
\tilde u(\xi)=A^{-1}_{\neq}(\xi)Q_N(\xi)B_MQ^{-1}_N(\xi)A^{-1}_=(\xi)\widetilde{(\ell v)}(\xi)+
\]
\[
+A^{-1}_{\neq}(\xi)V^{-1}_{\varphi}\left(\sum\limits_{l_1=
1}^{n_1}\sum\limits_{l_2=1}^{n_2}\dots\sum\limits_{l_n=1}^{n_n}\tilde c_{L}(\xi'_{K})\xi^{l_1-1}_{k_1}\xi^{l_2-1}_{k_2}\cdots\xi^{l_n-1}_{k_n}\right),
\]
where $c_L(x'_K)\in H^{S_L}({\mathbb R}^{M-n})$ are arbitrary functions,  $S_L=(s_1-\ae_1+l_1-1/2,\dots,s_n-\ae_n+l_n-1/2),~ l_j=1,2,\dots,n_j),~j=1,2,...,n,$ $\ell v$ is an arbitrary continuation of $v$ onto $H^{S-\alpha}({\mathbb R}^M)$.

The a priori estimate
\[
||u||_S\leq~const~\left(||v||^+_{S-\alpha}+\sum\limits_{l_1=
1}^{n_1}\sum\limits_{l_2=1}^{n_2}\dots\sum\limits_{l_n=1}^{n_n}||c_L||_{S_L}\right)
\]
holds.
}

\begin{proof}
Similar to the proof of Theorem 1 we obtain the equality (\ref{4}). Further,
let us note that the function  $A^{-1}_=(\xi)\widetilde{(\ell v)}(\xi)$  belongs to the space  $\tilde H^{S-\ae}({\mathbb R}^M)$. So, if take an arbitrary polynomial  $Q_N(\xi)$ satisfying the condition (\ref{5})
then the function $Q_N^{-1}(\xi)A^{-1}_=(\xi)\widetilde{(\ell v)}(\xi)$ will belong to the space  $\tilde H^{-\delta}({\mathbb  R}^M)$.

Further, according to the theory of multidimensional Riemann problem  \cite{V0} we can represent the latter function as a sum of two summands, this is so called a jump problem which can be solved by the operator  $B_M$:
$$
Q_N^{-1}A^{-1}_=\widetilde{(\ell v)}=f_++f_-,
$$
where $f_+\in\tilde H^{-\delta}(C), f_-\in\tilde H^{-\delta}({\mathbb R}^M\setminus C),$
\[
f_+=B_M(A_=^{-1}\widetilde{(\ell v)}),~~~f_-=(I-B_M)(A_=^{-1}\widetilde{(\ell v)}).
\]

Multiplying the equality  (\ref{4}) by $Q_N^{-1}(\xi)$ we rewrite it in the form
$$
Q_N^{-1}A_{\neq}\tilde u+Q_N^{-1}A_{=}^{-1}\tilde u_-=f_++f_-,
$$
or
$$
Q_N^{-1}A_{\neq}\tilde u-f_+=f_--Q_N^{-1}A_{=}^{-1}\tilde u_-
$$

In other words
\begin{equation}\label{6}
A_{\neq}\tilde u-Q_Nf_+=Q_Nf_--A_=^{-1}\tilde u_-.
\end{equation}

The left hand side of the equality  (\ref{6}) belongs to the space  $\tilde H^{-N-\delta}(C)$, bur the right hand side belongs to the space  $\tilde H^{-N-\delta}({\mathbb R}^M\setminus C)$. Therefore, we have
$$
F_M^{-1}(A_{\neq}\tilde u-Q_Nf_+)=F_M^{-1}(Q_Nf_--A_=^{-1}\tilde  u_-),
$$
where the left hand side belongs to the space $ H^{-N-\delta}(C)$, but right hand side belongs to the space $ H^{-N-\delta}({\mathbb R}^M\setminus C)$, from which we conclude immediately that this is a distribution supported on the surface $\partial C$.

The form for such a distribution is given in \cite{V6} for the cone $C_{K_j}$ with help of the operator $V_{\varphi_j}$. Thus,  we apply the operator $T_{\varphi}$ to the latter equality and obtain
$$
T_{\varphi}F_M^{-1}(A_{\neq}\tilde u-Q_Nf_+)=T_{\varphi}F_M^{-1}(Q_Nf_--A_=^{-1}\tilde  u_-),
$$
so that both left hand side and right hand side is a distribution supported on the hyper-plane $x_{k_1}=0, x_{k_2}=0,\dots,x_{k_n}=0$. Then
\[
T_{\varphi}F_M^{-1}(A_{\neq}\tilde u-Q_Nf_+)=
\]
\[
=\sum\limits_{l_1=
1}^{n_1}\sum\limits_{l_2=1}^{n_2}\dots\sum\limits_{l_n=1}^{n_n}c_{L}(x'_{K})\delta^{(l_1-1)}(x_{k_1})\delta^{(l_2-1)}(x_{k_2})\cdots\delta^{(l_n-1)}(x_{k_n}),
\]
where $L=l_1,\dots,l_n), x'_K=(x'_{K_1},\dots,x'_{K_n})\in{\mathbb R}^{M-n}, \delta$ is the Dirac mass-function.

Applying the Fourier transform we obtain
\[
F_MT_{\varphi}F_M^{-1}(A_{\neq}\tilde u-Q_Nf_+)=
\]
\begin{equation}\label{8}
=\sum\limits_{l_1=
1}^{n_1}\sum\limits_{l_2=1}^{n_2}\dots\sum\limits_{l_n=1}^{n_n}\tilde c_{L}(\xi'_{K})\xi^{l_1-1}_{k_1}\xi^{l_2-1}_{k_2}\cdots\xi^{l_n-1}_{k_n},
\end{equation}
Taking into account that $F_MT_{\varphi}F_M^{-1}$ we can write
\[
A_{\neq}\tilde u-Q_Nf_+
=V^{-1}_{\varphi}\left(\sum\limits_{l_1=
1}^{n_1}\sum\limits_{l_2=1}^{n_2}\dots\sum\limits_{l_n=1}^{n_n}\tilde c_{L}(\xi'_{K})\xi^{l_1-1}_{k_1}\xi^{l_2-1}_{k_2}\cdots\xi^{l_n-1}_{k_n}\right),
\]
or finally
\begin{equation}\label{10}
\begin{array}{rcl}
\tilde u(\xi)=A^{-1}_{\neq}(\xi)Q_N(\xi)B_MQ^{-1}_N(\xi)A^{-1}_=(\xi)\widetilde{(\ell v)}(\xi)+\\
+A^{-1}_{\neq}(\xi)V^{-1}_{\varphi}\left(\sum\limits_{l_1=
1}^{n_1}\sum\limits_{l_2=1}^{n_2}\dots\sum\limits_{l_n=1}^{n_n}\tilde c_{L}(\xi'_{K})\xi^{l_1-1}_{k_1}\xi^{l_2-1}_{k_2}\cdots\xi^{l_n-1}_{k_n}\right).
\end{array}
\end{equation}

To obtain a priori estimates let us note that all summands in the formula (\ref{8}) should belong to the space  $\widetilde H^{S-\ae}({\mathbb R}^M)$. We take one of summands and estimate corresponding integral.
\[
||\tilde c_{L}(\xi'_{K})\xi^{l_1-1}_{k_1}\xi^{l_2-1}_{k_2}\cdots\xi^{l_n-1}_{k_n}||^2_{S-\ae}\leq
\]
\[
\leq\int\limits_{{\mathbb R}^M}|\tilde c_{L}(\xi'_{K})|^2\prod\limits_{j=1}^n(1+|\xi_{K_j})^{2(s_j-\ae_j)}\prod\limits_{j=1}^n|\xi_{k_j}|^{2(l_j-1)}d\xi'_K\prod\limits_{j=1}^nd\xi_{k_j}\leq
\]
\[
\leq\int\limits_{{\mathbb R}^M}|\tilde c_{L}(\xi'_{K})|^2\prod\limits_{j=1}^n(1+|\xi_{K_j}|^{2(s_j-\ae_j+l_j-1)}d\xi'_K\prod\limits_{j=1}^nd\xi_{k_j},
\]
and for existence of each integral of the type
\[
\int\limits_{-\infty}^{+\infty}(1+|\xi'_{K_j}|+|\xi_{k_j}|)^{2(s_j-\ae_j+l_j-1)}d\xi_{k_j}
\]
the condition
\begin{equation}\label{9}
2(s_j-\ae_j+l_j-1)<-1
\end{equation}
is necessary. It is equivalent to the following condition
\[
s_j-\ae_j+l_j<1.
\]
Since we have $s_j-\ae_j+l_j=-n_j-\delta_j+l_j$  then we see that the condition (\ref{9}) is satisfied for all
\[
l_j=1,2,\dots,n_j,
\]
but it is not satisfied for $l_j=n_j+1$. After integration on all $\xi_{k_j}$ we will find that $\tilde c_L(\xi'_K)\in\widetilde H^{S_L}({\mathbb R}^{M-n})$, where $S_L=(s_1-\ae_1+l_1-1/2,\dots,s_n-\ae_n+l_n-1/2),$ and $l_j=1,2,\dots,n_j, j=1,2.\dots,n$.

For a priori estimates we have
\[
||A^{-1}_{\neq}(\xi)Q_N(\xi)B_MQ^{-1}_N(\xi)A^{-1}_=(\xi)\widetilde{(\ell v)}(\xi)||_S\leq
\]
\[
\leq~const~||B_MQ^{-1}_N(\xi)A^{-1}_=(\xi)\widetilde{(\ell v)}(\xi)||_{S-\ae+N}\leq
\]
\[
\leq~const~||Q^{-1}_N(\xi)A^{-1}_=(\xi)\widetilde{(\ell v)}(\xi)||_{S-\ae+N}\leq
\]
\[
\leq~const~||\widetilde{(\ell v)}(\xi)||_{S-\ae+N-N+\ae-\alpha}=const~||\widetilde{(\ell v)}(\xi)||_{S-\alpha}\leq~const~|| v||^+_{S-\alpha}
\]
according to Lemma 1 and the fact that $S-\ae+N=\delta. |\delta_j|<1/2, j=1,\dots,n$.

To estimate other summands in the formula (\ref{10}) we use above considerations. Really, if $\tilde c_L(\xi'_K)\in\widetilde H^{S_L}({\mathbb R}^{M-n})$ then each summand $\tilde c_{L}(\xi'_{K})\xi^{l_1-1}_{k_1}\xi^{l_2-1}_{k_2}\cdots\xi^{l_n-1}_{k_n}$ in the formula (\ref{10}) belongs to the space $\widetilde H^{S-\ae}({\mathbb R}^M)$. Thus, we have
\[
||A^{-1}_{\neq}(\xi)V^{-1}_{\varphi}\tilde c_{L}(\xi'_{K})\xi^{l_1-1}_{k_1}\xi^{l_2-1}_{k_2}\cdots\xi^{l_n-1}_{k_n}||_{S}\leq
\]
\[
\leq~const~||V^{-1}_{\varphi}\tilde c_{L}(\xi'_{K})\xi^{l_1-1}_{k_1}\xi^{l_2-1}_{k_2}\cdots\xi^{l_n-1}_{k_n}||_{S-\ae}\leq
\]
\[
\leq~const~||\tilde c_{L}(\xi'_{K})\xi^{l_1-1}_{k_1}\xi^{l_2-1}_{k_2}\cdots\xi^{l_n-1}_{k_n}||_{S-\ae}\leq~const~||\tilde c_L||_{S_L}.
\]
The latter estimate was obtained above. The Theorem 2 is proved.
\qed
\end{proof}

{\bf Remark 1.}
This formula includes the operator $V_{\varphi}$. Examples 2 and 3 give exact representation for this operator for certain concrete cones.

\section*{Conclusion}

These studies led to different boundary value problems for such
elliptic pseudo-differential equations in cones similar to \cite{V0,V1,V2}.
Particularly, for the case of Theorem 2 a general solution of the equation
(\ref{3}) includes a lot of arbitrary functions from corresponding Sobolev--Slobodetskii
spaces. To determine these functions uniquely one needs some additional conditions
(not necessary boundary conditions). We will try to describe certain
 statements of boundary value problems in forthcoming papers.
%


\end{document}